\documentclass[a4paper,11pt]{article}
\usepackage{array}
\usepackage{theorem}
\usepackage{amsmath,amscd,amssymb}
\usepackage{latexsym}
\usepackage{epsfig}

\theorembodyfont{\sl}

\newtheorem{lemma}{Lemma}[section]
\newtheorem{proposition}[lemma]{Proposition}
\newtheorem{theorem}[lemma]{Theorem}

\renewcommand{\AA}{\mathbb A}

\newcommand{\CC}{\mathbb C}

\newcommand{\PP}{\mathbb P}
\newcommand{\QQ}{\mathbb Q}
\newcommand{\RR}{\mathbb R}

\newcommand{\ZZ}{\mathbb Z}

\newcommand{\cD}{\mathcal D}

\newcommand{\cF}{\mathcal F}

\newcommand{\cL}{\mathcal L}
\newcommand{\cM}{\mathcal M}

\newcommand{\cO}{\mathcal O}

\renewcommand{\Tilde}{\widetilde}

\renewcommand{\Bar}{\overline}

\newcommand{\GL}{\mathop{\mathrm {GL}}\nolimits}

\newcommand{\SO}{\mathop{\mathrm {SO}}\nolimits}
\newcommand{\Sp}{\mathop{\mathrm {Sp}}\nolimits}
\newcommand{\Orth}{\mathop{\null\mathrm {O}}\nolimits}
\newcommand{\POrth}{\mathop{\null\mathrm {PO}}\nolimits}
\newcommand{\POrthtd}{\mathop{\null\mathrm {P\widetilde{O}}}\nolimits}
\newcommand{\PGamma}{\mathop{\null\mathrm {P\Gamma}}\nolimits}
\newcommand{\Proj}{\mathop{\null\mathrm {Proj}}\nolimits}

\newcommand{\Aut}{\mathop{\mathrm {Aut}}\nolimits}

\newcommand{\Hom}{\mathop{\mathrm {Hom}}\nolimits}
\newcommand{\Ker}{\mathop{\mathrm {Ker}}\nolimits}

\newcommand{\Mat}{\mathop{\mathrm {Mat}}\nolimits}

\newcommand{\rank}{\mathop{\mathrm {rank}}\nolimits}
\newcommand{\norm}{\mathop{\mathrm {norm}}\nolimits}
\newcommand{\scale}{\mathop{\mathrm {scale}}\nolimits}

\newcommand{\even}{\mathop{\mathrm {even}}\nolimits}
\newcommand{\odd}{\mathop{\mathrm {odd}}\nolimits}
\newcommand{\tor}{\mathop{\mathrm {tor}}\nolimits}

\newcommand{\vol}{\mathop{\mathrm {vol}}\nolimits}
\newcommand{\geom}{\mathop{\mathrm {geom}}\nolimits}

\newcommand{\trace}{\mathop{\mathrm {tr}}\nolimits}

\newcommand{\id}{\mathop{\mathrm {id}}\nolimits}
\newcommand{\sn}{\mathop{\mathrm {sn}}\nolimits}

\newcommand{\bdv}{\mathbf v}
\newcommand{\bdw}{\mathbf w}
\newcommand{\bdx}{\mathbf x}
\newcommand{\bdy}{\mathbf y}

\newcommand{\latt}[1]{{\langle{#1}\rangle}}

\newcommand{\Kthree}{\mathop{\mathrm {K3}}\nolimits}

\newcommand{\Lrank}{\rho}
\newcommand{\qedsymbol}{\mbox{$\Box$}}
\newcommand{\qed}{\unskip\nobreak\hfil\penalty50\hskip1em\hbox{}\nobreak
\hfill\qedsymbol\parfillskip=0pt\finalhyphendemerits=0}
\newenvironment{proof}{\begin{ProofwCaption}{Proof}}{\end{ProofwCaption}}
\newenvironment{ProofwCaption}[1]
 {\addvspace\theorempreskipamount \noindent{\it #1.}\rm}
 {\qed \par \addvspace\theorempostskipamount}

\begin{document}

\title{The Hirzebruch-Mumford volume for the orthogonal group
and applications}
\author{V. Gritsenko, K.~Hulek and G.K. Sankaran}
\date{}
\maketitle

\begin{abstract}
In this paper we derive an explicit formula for the Hirzebruch-Mumford
volume of an indefinite  lattice $L$ of rank $\ge 3$.
If $\Gamma  \subset \Orth(L)$ is an arithmetic subgroup
and $L$ has signature $(2,n)$,
then an application of  Hirzebruch-Mumford proportionality allows us
 to determine the leading
term of the growth of the dimension of the spaces $S_k(\Gamma)$ of
cusp forms of weight $k$, as $k$ goes to infinity. We compute this in
a number of examples, which are important for geometric applications.
\end{abstract}

\section{Introduction}\label{Intro}

It is well known that one can use the Hirzebruch-Mumford
proportionality principle (see Theorem \ref{hmumford1} and Corollary
\ref{hmumford2}) in order to estimate the growth behaviour, as the
weight increases, of the dimension of spaces of modular forms for some
arithmetic group $\Gamma$ acting on a homogeneous domain. Often, a
qualitative statement suffices, but in other cases one needs the exact
form of the leading coefficient of the Hilbert polynomial. For this,
a precise knowledge of the constant of proportionality that appears
in the Hirzebruch-Mumford theorem is necessary. This in turn
requires the computation of what we call the {\em Hirzebruch-Mumford
volume} of the arithmetic group $\Gamma$. If $\Gamma$ acts freely,
this can be defined as the quotient of the Euler number of a
fundamental domain of $\Gamma$ by the Euler number of the compact
dual.

The subject starts with the seminal work of Siegel \cite{Sie1}
on the volume of
the orthogonal group. Very many authors have taken up his theory and
generalised it in many different directions, including Harder
\cite{Ha}, Serre \cite{Se}, Prasad \cite{Pr} and many others. Our
specific interest lies in indefinite orthogonal groups (see the work
by Shimura \cite{Sh}, Gross \cite{Gr}, Gan, Hanke and Yu \cite{GHY},
as well as Belolipetsky and Gan \cite{BG}, to name some important
recent work in this direction). Motivated by possible applications
(cf. \cite{GHS1}, \cite{GHS2}) concerning moduli spaces of $\Kthree$ surfaces and
similar modular varieties we started to investigate the volume of
certain arithmetic subgroups of orthogonal groups $\Orth(L)$ of even
indefinite lattices of signature $(2,n)$. All our groups are defined
over the rational numbers; but for the applications we have in mind
we cannot restrict ourselves to
unimodular or maximal lattices. To our knowledge there exist no
results in the literature that allow an easy calculation of the
Hirzebruch-Mumford volume for the groups we treat in this paper.

In order to compute these volumes we therefore decided to return to
Siegel's work. Let $L$ be an even indefinite lattice of signature
$(2,n)$ and let $\Orth(L)$ be its group of isometries. The lattice $L$
defines a domain
$$
\Omega_{L}=\{\left[\bdw \right]\in\PP(L\otimes \CC); \;
(\bdw,\bdw)_{L}=0, \;
(\bdw,\Bar{\bdw} )_{L}>0\}.
$$
This domain has two connected components $\cD_{L}$ and $\cD'_{L}$,
which are interchanged by complex conjugation. Let $\Orth^+(L)$ be
the index $2$ subgroup of $\Orth(L)$ which fixes $\cD_L$. The
fundamental problem of our paper is to determine the
Hirzebruch-Mumford volume of this group. To do this, one has to
compare the volume of the quotient $\Orth^+(L) \backslash \cD_L$ to
the volume of the compact dual $\cD_L^{(c)}$. To do so correctly, one
has to choose volume forms on the domain $\cD_L$ and the compact dual
$\cD_L^{(c)}$ that coincide at the common point of both domains given
by a maximal compact subgroup. Note that here
$\cD_L=\Orth(2,n)/\Orth(2) \times \Orth(n)$ and $\cD_L^{(c)}=
\Orth(2+n)/\Orth(2) \times \Orth(n)$. This is in fact a problem which
does not depend on the complex structure of the domains, but can be
considered in greater generality for indefinite lattices of signature
$(r,s)$. We use the volume form on $\cD_L$ which was introduced by
Siegel. It then turns out that this must be compared to the volume
form on $\cD_L^{(c)}$ which is given by $1/2$ of the volume form
induced by the Killing form on the Lie algebra of the group
$\SO(r+s)$. Comparing these two volumes gives us the main formula for
the Hirzebruch-Mumford volume of $\Orth^+(L)$ This formula involves
the Tamagawa (Haar) measure of the group $\Orth(L)$. However, again
using a result of Siegel, the computation of the Tamagawa measure can
be reduced to computing the local densities $\alpha_p(L)$ of the lattice $L$
over  the $p$-adic integers. Our main formula
for any indefinite lattice $L$ of rank $\Lrank\ge 3$ is
$$
\vol_{HM}(\Orth(L))=\frac{2}{g_{sp}^+}\,|\det L|^{(\Lrank+1)/2}
\prod_{k=1}^{\Lrank}\pi^{-k/2}\Gamma(k/2)
 \prod_{p}\alpha_p(L)^{-1}
$$
where $g_{sp}^+$  is the number of the proper spinor genera in the genus
of $L$ (see Theorem  \ref{prop_mainformula}).
Since everything is defined over the
rationals, one can use Kitaoka's book \cite{Ki} on quadratic forms to
compute the local densities in question.

In order to illustrate our results, and particularly in view of
applications, we compute the Hirzebruch-Mumford volume for several
examples. The lattices and the groups which we consider are mostly
related to moduli problems. We start with a series of even unimodular
lattices, namely the lattices ${II}_{2,2m+8}=2U \oplus mE_8(-1)$, where $U$
denotes the hyperbolic plane and $E_8$ is the positive definite root
lattice associated to $E_8$. The next series of examples consists of
the lattices $L_{2d}^{(m)}=2U\oplus mE_8(-1)\oplus \latt{-2d}$, which
are closely related to well known moduli problems. Let
$$
{\cF}_{2d}^{(m)} = \widetilde{\Orth}^+(L_{2d}^{(m)}) \backslash \cD_{L_{2d}^{(m)}}
$$
where $\widetilde{\Orth}^+(L_{2d}^{(m)})$ is the subgroup of
$\Orth^+(L_{2d}^{(m)})$ which acts trivially on the discriminant
group.
For $m=0$ and $d$ a prime number, ${\cF}_{2d}^{(0)}$ is a moduli of Kummer
surfaces  (see \cite{GH}). The spaces ${\cF}_{2d}^{(1)}$ parametrise certain
lattice-polarised  $\Kthree$ surfaces and
if $m=2$,  then ${\cF}_{2d}={\cF}_{2d}^{(2)}$ is the moduli space of
$\Kthree$ surfaces of degree $2d$.
We compute the volumes of the groups $\Orth^+(L_{2d}^{(m)})$
and $\widetilde{\Orth}^+(L_{2d}^{(m)})$ and obtain as a corollary the
leading term controlling the growth behaviour of the dimension of the
spaces of cusp forms for these groups. As a specialisation of this
example we recover known formulae for the Siegel modular group in
genus $2$ and the paramodular group.
The series of examples considered in this paper, namely
the even indefinite unimodular lattices (Section \ref{evenunimodular}), their sublattices
$T$ (Section \ref{Tm}) and some lattices of signature $(2,8m+2)$ (Section \ref{Heegner}),
are closely related to moduli of $\Kthree$ surfaces and related quotients of
homogeneous varieties of type IV.
The volumes of these lattices determine the part of the obstruction for extending
pluricanonical differential forms on ${\cF}_{2d}^{(m)}$
to a smooth compactification of this variety which comes from the ramification divisor.

In~\cite{GHS2} we use these results to obtain information about the Kodaira
dimension of two series of modular varieties, including effective bounds on the degree $d$
which guarantee that the varieties ${\cF}_{2d}^{(m)}$ are of general
type. The case of polarised $\Kthree$ surfaces is considered in~\cite{GHS1}.
\footnote{For more details of this, see the talk at
the Arbeitstagung in Bonn by
the first author on 13th June 2005: {\tt
http://www.mpim-bonn.mpg.de/at2005}.}

The paper is organised as follows: in Section
\ref{sec_HMproportionality} we recall Hirzebruch-Mumford
proportionality and the Hirzebruch-Mumford volume in the form in which
we need it (see Theorem \ref{hmumford1} and Corollary
\ref{hmumford2}). In Section \ref{Computationvolumes} we perform the
necessary volume computations and derive the main formula (see Theorem
\ref{prop_mainformula}). In Section \ref{Applications} we treat in
some detail several lattices which appear naturally in moduli
problems.

\subsection*{Acknowledgements}

{\small We are grateful to various organisations that have supported
us during the preparation of this work. The DFG Schwerpunktprogramm
SPP 1084 ``Globale Me\-tho\-den in der komplexen Geometrie'', grant HU
337/5-2, has enabled VAG and GKS to visit the University of Hannover
for extended stays. VAG would like to thank the Max-Planck-Institut
f\"ur Mathematik in Bonn for hospitality in 2005. KH is grateful to
the Fields Institute in Toronto and the Graduate School of Mathematics
of Nagoya University, in particular to Professors N.~Yui and S.~Kondo.
GKS is grateful to Tokyo University and to the Royal Society. All
these institutions provided excellent working conditions. We should
like to thank M.~Belolipetsky and R.~Schulze-Pillot for useful discussions
on the mass formula.}

\section{Hirzebruch-Mumford proportionality}\label{sec_HMproportionality}

In this section we consider an indefinite even lattice $L$ of
signature $(2,n)$. Let $\Orth(L)$ be its group of isometries. We
denote by $(\ ,\ )_{L}$ the form defined on $L$, extended bilinearly
to $L \otimes \RR$ and $L \otimes \CC$. The domain
$$
\Omega_{L}=\{\left[\bdw \right]\in\PP(L\otimes \CC); \;
(\bdw,\bdw)_{L}=0, \;
(\bdw,\Bar{\bdw} )_{L}>0\}
$$
has two connected components, say $\Omega_{L}=\cD_{L}\cup
\cD'_{L}$, which are interchanged by complex conjugation. By
$\cD_{L}^{\bullet}$ we denote the affine cone over $\cD_{L}$ in $L
\otimes \CC$. Let $\Gamma \subset \Orth(L)$ be an arithmetic group
which leaves the domain $\cD_{L}$ invariant. A {\em modular form\/} of
{\em weight} $k$ with respect to the group $\Gamma$ and with a (finite order)
{\em character} $\chi: \Gamma \to \CC^*$  is a holomorphic map
$$
f: \cD_{L}^{\bullet} \to \CC
$$
which has the two properties
$$
\begin{array}{lcll}
f(tz)&=&t^{-k}f(z)& \mbox{for } t\in \CC^*,\\
f(g(z)) & = &\chi(g) f(z) & \mbox{for } g \in \Gamma. 
\end{array}
$$
If $n\le 2$ the function $f(z)$ must also be required to be holomorphic at
infinity. A {\em cusp
form\/} is a modular form which vanishes on the boundary.

We denote the spaces of modular forms and of cusp forms of weight $k$,
with respect to the group $\Gamma$ and character $\chi$, by $M_k(\Gamma, \chi)$ and $S_k(\Gamma, \chi)$ 
respectively. These are finite dimensional vector spaces. 
Note that if $- \id \in \Gamma$ and $(-1)^k \neq \chi(- \id)$ then obviously 
$M_k(\Gamma, \chi)=0$. 
 
Modular forms can be interpreted as sections of
suitable line bundles. For this, we first assume that the group
$\Gamma$ is neat, in which case it acts freely on $\cD_L$ and we also assume that the character $\chi$
is trivial. Then the
transformation rules of modular forms of weight $1$ define a line
bundle $\cL$ on the quotient $\Gamma \backslash \cD_L$ and modular 
forms of weight $k$ with trivial character become sections in
$\cL^{\otimes k}$. The line 
bundle $\cL$, and its sections, extend to the Baily-Borel
compactification $\overline{\Gamma \backslash \cD}_L$. In fact, the
Baily-Borel compactification is the normal projective variety associated to
$\Proj\big(\bigoplus_k H^0(\cL^{\otimes k})\big)$. In general, modular forms of weight $k$ and with 
a character $\chi$ define sections of a line bundle $\cL_{k,\chi}$ which differs from
$\cL^{\otimes k}$ only by torsion. 

Every toroidal compactification $(\Gamma \backslash \cD_L)^{\tor}$
has a morphism $(\Gamma \backslash \cD_L)^{\tor} \to
\overline{\Gamma \backslash \cD}_L$ which is the identity on $\Gamma
\backslash \cD_L$. Via this morphism, we shall also consider $\cL$ and $\cL_{k,\chi}$
as line bundles on $(\Gamma \backslash \cD_L)^{\tor}$, using the same
symbol by abuse of notation. If $\Gamma$ is not neat then the above
remains true, as long as we consider $\cL$ and $\cL_{k,\chi}$ as $\QQ$-line
bundles or only consider weights $k$ that are sufficiently divisible.

The connection with pluricanonical forms is as follows.  
There is an $n$-form $dZ$ on $\cD_L$ such that if $f$ is a modular
form of weight $n = \dim \cD_L$ with character $\det$, then $\omega=f dZ$ is a 
$\Gamma$-invariant $n$-form on $\cD_L$. Hence, if the action of
$\Gamma$ on $\cD_L$ is free, $\omega$ descends to an $n$-form on
$\Gamma \backslash \cD_L$. Similarly, modular forms of weight $kn$ with character $\det^k$ 
define $k$-fold pluricanonical forms on $\Gamma \backslash \cD_L$. If
$\Gamma$ does not act freely, then this is still true outside the
ramification locus of the quotient map $\cD_L \to \Gamma \backslash
\cD_L$. These forms will, in general, not extend to compactifications
of $\Gamma \backslash \cD_L$. If $\Gamma$ is a neat group, then let
$(\Gamma \backslash \cD_L)^{\tor}$ be a smooth toroidal
compactification (which always exists by \cite{AMRT}). Let $D$ be the
boundary of such a toroidal compactification. If $\cM$ is the line bundle
of modular forms of weight $n$ and character $\det$,  
then the canonical bundle is 
given by $\omega_{(\Gamma \backslash \cD_L)^{\tor}} = \cM
\otimes \cO_{(\Gamma \backslash \cD_L)^{\tor}}(-D)$. Hence, if $f$ is
a weight $n$ form with character $\det$, not vanishing at the boundary, then $f dZ$ defines 
an $n$-form on $(\Gamma \backslash \cD_L)^{\tor}$ with poles along the
boundary. However, if $f$ is a cusp form, then $f dZ$ does define an
$n$-form on $(\Gamma \backslash \cD_L)^{\tor}$, and similarly forms of
weight $kn$ and character $\det^k$, which vanish along the boundary of order $k$, define 
$k$-fold pluricanonical forms on $(\Gamma \backslash
\cD_L)^{\tor}$. It should be pointed out that some authors define
automorphic forms a priori as those functions that give rise to
pluricanonical forms. In our context, this means a restriction to
forms of weight $kn$. Moreover, the weight of these forms is sometimes
defined as $k$. We shall refer to the latter as the {\em geometric\/}
weight, in contrast to the {\em arithmetic\/} weight of our
definition. This difference accounts for the fact that some of our
formulae differ from corresponding formulae in the literature by
powers of $n$.

The Hirzebruch-Mumford proportionality principle, which works very
generally for quotients of a homogeneous domain $\cD$ by an arithmetic
group $\Gamma$, allows us to estimate the growth behaviour of spaces
of cusp forms as a function of the weight $k$ in terms of a suitably
defined volume. This was first discovered by Hirzebruch \cite{Hi1},
\cite{Hi2} in the case where the quotient $\Gamma \backslash \cD$ is
compact, and was generalised by Mumford \cite{Mum} to the case where
$\Gamma \backslash \cD$ has finite volume. We denote the compact dual
of $\cD$ by $\cD^{(c)}$. Let $\overline{X}$ be the Baily-Borel
compactification of $X=\Gamma \backslash \cD$ and let $X^{\tor}$ be
some smooth toroidal compactification of $X$.

\begin{theorem}\label{hmumford1}
Let $\Gamma$ be a neat arithmetic group which acts on a bounded
symmetric domain $\cD$.
Let $S^{\geom}_k (\Gamma)=S_{nk} (\Gamma, \det^k)$ be the space of 
cusp forms of geometric weight $k$ with respect to $\Gamma$. Then
$$
\dim S^{\geom}_k (\Gamma)=
\vol_{HM}(\Gamma)h^0(\omega^{(1-k)}_{\cD^{(c)}})+P_1(k)
$$
where $P_1(k)$ is a polynomial whose degree is at most the dimension
of $\Bar X \backslash X$.
\end{theorem}
\begin{proof}
This is \cite[Corollary 3.5]{Mum}.
\end{proof}

Here $\vol_{HM}(\Gamma \backslash \cD)$ denotes a suitably normalised
volume of the quotient $\Gamma \backslash \cD$, which we shall refer
to as the {\em Hirzebruch-Mumford volume}. If $\Gamma$ acts freely,
then the Hirzebruch-Mumford volume is a quotient of Euler numbers
$$
\vol_{HM}(\Gamma)=\vol_{HM}(\Gamma \backslash \cD)= \frac{e(\Gamma
 \backslash \cD)}{e(\cD^{(c)})}.
$$
If $\Gamma$ does not act freely, then choose a normal subgroup
$\Gamma' \lhd \Gamma$ of finite index which does act freely. Then
$$
\vol_{HM}(\Gamma)= \frac{\vol_{HM}(\Gamma')}{[\PGamma:\Gamma']}
$$
where $\PGamma$ is the image of $\Gamma$ in $\Aut(\cD)$, i.e. the
group $\Gamma$ modulo its centre. This value is independent of the
choice of the subgroup $\Gamma'$.

We would like to point out that the asymptotic growth of the dimension of the space of modular forms
is the same as that for cusp forms. This follows since the difference can be bounded by the dimension of
a space of sections of a line bundle supported on the boundary.

Hirzebruch \cite{Hi1} first formulated his result in the case where the group is cocompact, i.e.,
where the quotient
$X=\Gamma \backslash \cD$ is compact. Since the Chern numbers of $X$ and that of the compact dual
are proportional and the factor of proportionality is given by the volume, one can use Riemann-Roch to
compute the
exact dimension of the space of modular foms (in this case it does not make sense to talk about cusp forms).

We shall now apply this to orthogonal lattices.

\begin{proposition}\label{hmumford2}
Let $L$ be an indefinite even lattice of signature $(2,n)$ and let
$\Gamma$ be an arithmetic subgroup which acts on the domain $\cD$.
Fix a positive integer $k$ and a character $\chi$. If $- \id \in \Gamma$, then we 
restrict to those $k$ for which $(-1)^k= \chi(- \id)$. Then the dimension of 
the space $S_k(\Gamma, \chi)$ of cusp forms of arithmetic weight $k$
grows 
as
$$
\dim S_k(\Gamma, \chi)=\frac{2}{n!}\vol_{HM}
(\Gamma\backslash\cD_L)k^n+O(k^{n-1}).
$$
\end{proposition}
\begin{proof}
We shall first assume that $\Gamma$ is neat (in which case
automatically $- \id \notin \Gamma$) and that $\chi$ is trivial.  
We consider $\cL$ as a 
line bundle on a smooth toroidal compactification $X^{\tor}$ of
$X=\Gamma\backslash\cD_L$. 
It follows from the definition of cusp forms that $H^0 (X^{\tor},{\cL}^{\otimes k}(-D))=S_k(\Gamma)$. 
Since $\cL$ is big and nef and $K_{X^{\tor}}={\cL}^{\otimes n}(-D)$, 
it follows from Kawamata-Viehweg vanishing 
that $h^i(X^{\tor}, \cL^{\otimes k}(-D))=0$ for $i\geq 1$ and $k\gg 0$ and
hence $\chi(X^{\tor}, \cL^{\otimes k}(-D))=h^0(X^{\tor}, \cL^{\otimes k}(-D))$ for $k\gg 0$. 
The leading term of the Riemann-Roch polynomial as a function of $k$ is given by
$c_1^n(\cL)/n!$. The same argument goes through for $\cL_{k,\chi}$.
Since $\cL^{\otimes k}$ and $\cL_{k,\chi}$ only differ by torsion they have the 
same leading coefficients.  

In order to apply Theorem \ref{hmumford1} we consider the line bundle
$\cM$, whose sections are modular forms of geometric 
weight $1$, i.e modular forms of arithmetic weight $n$ and character $\det$. 
Note that $\cM^k=\cL^{nk}$ for suitably divisible $k$. 
Also recall that in the orthogonal case the compact dual 
$\cD^{(c)}$ is the $n$-dimensional quadric. Hence, by a
straightforward calculation, the leading term of
$h^0(\omega^{(1-k)}_{\cD^{(c)}})$ can be seen to be equal to
$2n^n/n!$. It then follows from Hirzebruch-Mumford proportionality
that
$$
\frac{c^n_1 (\cM^n)}{n!}= \frac{2n^n}{n!} \vol_{HM}(X)
$$
and hence
$$
\frac{c^n_1(\cL)}{n!}=\frac{2}{n!}\vol_{HM}(X)
$$
which gives the claim in the case of a neat group.

We now consider a group $\Gamma$ which is not necessarily neat and
choose $\Gamma' \lhd \Gamma$ neat and of finite index. The group
$\Gamma$ acts on the total space of the line bundle $\cL$, and if $-
\id \in \Gamma$ then it follows from our assumptions on $k$  
that this element acts trivially. We can now apply the Lefschetz fixed
point formula (cf. \cite[Appendix to \S 2]{T}), from which we obtain
\begin{eqnarray*}
\dim S_k (\Gamma)&=&\dim S_k(\Gamma')^{\Gamma}\\
&=&\frac{1}{[\PGamma:\Gamma']}\cdot\sum\limits_{\gamma\in\PGamma/\Gamma'}
\trace\left(\gamma|_{S_k(\Gamma')}\right)\\
&=&\frac{1}{[\PGamma:\Gamma']}\dim S_k(\Gamma')+O(k^{n-1})\\
&=&\frac{1}{[\PGamma:\Gamma']}\vol_{HM}(\Gamma'\backslash\cD_L)
\frac{2}{n!}k^n+O(k^{n-1})\\
&=& \frac{2}{n!}\vol_{HM}(\Gamma\backslash\cD_L)k^n+O(k^{n-1}).
\end{eqnarray*}
\end{proof}

Note that the growth behaviour of the space of modular forms of weight $k$ and that of the 
space of cusp forms are the same. This follows from the exact sequence
$$
0\to \cL^{\otimes k} (-D)\to \cL^{\otimes k} \to \cL^{\otimes k}|_D \to 0.
$$

\section{Computation of volumes}\label{Computationvolumes}

In order to compute the leading coefficient which determines the
growth of the dimension of spaces of cusp forms, we have to compare
the volume of a fundamental domain of an arithmetic group $\Gamma$ to
the volume of the compact dual. For this, the complex structure is not
important and we therefore consider, more generally, an indefinite
integral lattice $L$ of signature $(r,s)$.

As before, we denote the group of isometries of the lattice $L$ by
$\Orth(L)$. The lattice $L$ defines  a homogeneous domain
$\cD_{rs}$.
In terms of groups the domain
$\cD_{rs}$ is the quotient of the orthogonal group $\Orth(L \otimes
\RR)$ by a maximal compact subgroup, i.e.,
$$
\cD_{rs}=\cD_L=\Orth(r,s)/\Orth(r) \times \Orth(s)=\SO^+(r,s)/\SO(r)
\times \SO(s)
$$
where all groups are real Lie groups and $\SO^+(r,s)$ is the
connected component of the identity of $\SO(r,s)$.

The domain $\cD_{rs}$ can be realised as a bounded domain in the form
$$
\cD_{rs}=\{X\in \Mat_{r\times s}(\RR)\,;\, I_r-X{}^tX>0\}
$$
where $I_r\in \Mat_{r\times r}(\RR)$ is the identity matrix and the
action of the orthogonal group is given in the usual form, namely by
$$
M(X)=(AX+B)(CX+D)^{-1}
$$
for
$$
M=\begin{pmatrix}
A&B\\ C&D
\end{pmatrix} \in \Orth(r,s), \qquad A\in \Mat_{r\times r}(\RR),\
D\in\Mat_{s\times s}(\RR).
$$
We consider the $\Orth(r,s)$-invariant metric given by
$$
ds^2={\trace}\bigl((I_r-X{}^tX)^{-1}\,dX\, (I_s-{}^tXX)^{-1}\,
d{}^tX \bigr).
$$
Since
$$
\det((I_r-X{}^tX)^{-1})^s\cdot\det((I_s-{}^tXX)^{-1})^r=
\det((I_r-X{}^tX)^{-1})^{r+s}
$$
the corresponding volume form is given by
$$
dV = (\det(I_r-X{}^tX)^{-1})^{\frac{r+s}2}\,\prod_{i,j}dx_{ij}.
$$
Siegel computed the volume of $\cD_{rs}$ with respect to this volume
form in \cite{Sie2} (see also \cite[Theorem 7, p. 155]{Sie3}).
His result
is
\begin{equation}\label{eqn_Siegelvolume}
\vol_{S}(\Orth(L))
=\vol_{S}(\Orth(L) \backslash \cD_{rs})
=2\alpha_{\infty}(L)|\det L|^{(r+s+1)/2}\,
\gamma_r^{-1}\gamma_s^{-1}\,,
\end{equation}
where
\begin{equation}
\label{eqn_gamma}
\gamma_m=\prod_{k=1}^{m}\pi^{k/2}\,\Gamma(k/2)^{-1}
\end{equation}
and $\alpha_{\infty}(L)$ is the real Tamagawa (Haar) measure of the
lattice $L$. Formula~(\ref{eqn_Siegelvolume}) is valid for any
indefinite lattice
$L$ of rank $\ge 3$.
As indicated by the subscript, we shall refer to this
volume as the {\em Siegel volume\/} of the group $\Orth(L)$.

We want to understand the Siegel metric in terms of Lie algebras. Let
$\mathfrak g$ and $\mathfrak t$ be the Lie algebras of the indefinite
orthogonal group $\Orth(r,s)$ and its maximal compact subgroup
$\Orth(r)\times \Orth(s)$ respectively. Then
$$
\mathfrak g=\mathfrak t \oplus \mathfrak p
$$
where $\mathfrak p$ is the orthogonal complement of $\mathfrak t$
with respect to the Killing form. By \cite[p. 239]{He} this is
isomorphic to
$$
\mathfrak p=
\left\{\begin{pmatrix} 0&U\\{}^tU&0
\end{pmatrix}; \quad U\in \Mat_{r\times s}(\RR)\right\}.
$$
The space $\mathfrak p$ is isomorphic to the tangent space of
$\cD_{rs}$ at $0$. A straightforward calculation shows that the
$\Orth(r,s)$-invariant metric $ds^2$ is induced by the Killing
functional $\trace(U_1{}^tU_2)$ giving the scalar product
$\trace\bigl(dX_1{}^tdX_2\bigl)$ on the tangent space at $0$.

We now want to compare this to a suitable volume form on the compact
dual. Recall that the general situation is as follows. Let $H$ be a
bounded homogeneous domain and $G=\Aut(H)_0$ be the connected
component of the group of automorphisms of $H$. In particular,
$H=G/K$ where $K=G_{z_0}$ is the stabiliser of some point $z_0$.
There exists a unique compact real form $G_u$ of the complex group
$G_{\CC}$ such that $G\cap G_u =K$ and the symmetric domain $H=G/K$
can be embedded into the compact manifold $\cD^{(c)}=G_u/K$ as an open
submanifold. In our situation
$$
\cD^{(c)}_{rs}=\SO(r+s)/\SO(r)\times \SO(s).
$$
Again by \cite[p.239]{He} the tangent space of $\cD^{(c)}_{rs}$
at the point $I_{r+s}$ is given by the subspace
$$
\mathfrak p'=\left\{\begin{pmatrix} 0&U\\-{}^tU&0
\end{pmatrix}; \quad U\in \Mat_{r\times s}(\RR)\right\}
$$
of the Lie algebra of $\SO(r+s)$. The Killing form
$\trace(W_1{}^tW_2)$ of the Lie algebra of the compact group
$\SO(r+s)$ induces the form $2\trace(U_1 {}^tU_2)$ on the tangent
space $\mathfrak p'$. In order to compare the volumes of $\cD_{rs}$
and its compact dual $\cD^{(c)}_{rs}$ we have to normalise this form
in such a way that it coincides with the Siegel metric in the common
base point $K \in \cD_{rs} \subset \cD^{(c)}_{rs}$, i.e. we have to
use the form $\frac{1}{2}\trace(W_1{}^tW_2)$. Since the dimension of
$\SO(n)$ is $\frac{1}{2} {n(n-1)}$, we get a factor
$2^{-(r+s)(r+s-1)/4}$ in front of the volume of the compact group,
calculated in terms of the volume form induced by the Killing
functional on $\SO(r+s)$. The latter volume is computed in \cite[\S
3.7]{Hua}. Taking the above normalisation into account we find
\begin{equation}\label{eqn_volume(SO(m)}
\vol_S(\SO(m))=2^{m-1}\,\gamma_m
\end{equation}
and we shall again refer to this volume as the Siegel volume. For the
compact dual this gives
\begin{equation}\label{eqn_volumecompactdual}
\vol_S(\cD^{(c)}_{rs})=2\,\gamma_{r+s}\gamma_r^{-1}\gamma_s^{-1}.
\end{equation}
Our aim is to compute the Hirzebruch-Mumford volume
\begin{equation}\label{eqn_HMvolume}
\vol_{HM}(\Orth(L))=
\frac{\vol_{S}(\Orth(L)\backslash \cD_{rs})}{\vol_{S}(\cD^{(c)}_{rs})}.
\end{equation}
To make the above equation effective, we have to determine the
Tamagawa measure
$$
\alpha_{\infty}(L)=\alpha_{\infty}(\Orth(L)\backslash \Orth(L\otimes \RR))=
\alpha_{\infty}(\SO(L)\backslash \SO(L\otimes \RR)).
$$
The genus of the indefinite lattice $L$  contains
a finite number $g_{sp}^+(L)$ of (proper) spinor genera (for a definition see \cite[\S 6.3]{Ki}).
(We consider only proper classes and proper spinor genera.)
This number is always a power of two and can be calculated effectively.
It is well known that the spinor genus of an indefinite lattice of rank $\ge 3$
coincides with the class.
As was proved by
M. Kneser (see \cite{Kn})  the weight of the representations of
a given number $m$ by a spinor genus is the same for all genera
in the genus of $L$. The same arguments show that all spinor
genera in the genus
have the same mass.
(We are grateful to R.~Schulze-Pillot for drawing our attention to this fact.)
It is easy to see this in adelic terms. A spinor genus corresponds to
a double class
$\SO(V)\SO'_\AA(V)\, b\,\SO_\AA(L)$  in the adelic
group $\SO_\AA(V)$, where $V=L\otimes \QQ$ is the
rational quadratic space, and
$$
\SO'_\AA(V)=\ker \sn\colon \SO_\AA(V)\to
\QQ^\times_\AA/(\QQ^\times_\AA)^2
$$
is the kernel of the spinor norm. We note that
 the genus of $L$ is given by  $\SO_\AA(V)L$.
It follows from the definition that the group $\SO'_\AA(V)$
contains the commutator of $\SO_\AA(V)$, therefore
$$
\SO(V)\SO'_\AA(V)\, b\, \SO_\AA(L)= \SO(V)\SO'_\AA(V)\SO_\AA(L)\,b.
$$
The mass of a spinor genus
$$
\tau(\SO(V)\setminus \SO(V)\SO'_\AA(V) b \SO_\AA(L))=
\tau(\SO(V)\setminus \SO(V)\SO'_\AA(V) \SO_\AA(L))
$$
depends only on the genus, since the Tamagawa measure is invariant.
The  Tamagawa number of the orthogonal group is $2$
(see \cite{Sie1}, \cite{W}, \cite{Sh}), i.e.,
$\tau(\SO(V)\setminus \SO_\AA(V))=2$.
Then the Tamagawa measure $\alpha_{\infty}(L)$ can be
computed via the local densities of the lattices $L \otimes \ZZ_p$
over the $p$-adic integers $\ZZ_p$ (the local Tamagawa measures).
 More precisely,
\begin{equation}\label{eqn_localglobal} \alpha_\infty(L)=\alpha_\infty(\SO(L)\setminus \SO(L\otimes \RR))=
\frac{2}{g_{sp}^+(L)}\prod_p \alpha_p(L)^{-1},
\end{equation}
where $p$ runs through all prime numbers and $g_{sp}^+(L)$ is
the number of spinor genera in the genus of $L$.
The local densities can be computed, at least for
quadratic forms over $\QQ$ and its quadratic extensions: see
\cite{Ki}. In order to find $\alpha_p(L)$ it is enough to know the
Jordan decomposition of $L$ over the $p$-adic integers.

We can now summarise our results as follows
\begin{theorem}[Main formula]\label{prop_mainformula}
Let $L$ be an indefinite lattice of rank $\Lrank\ge 3$. Then the Hirzebruch-Mumford volume of $\Orth(L)$ equals
\begin{equation}\label{eqn_mainvolume}
\vol_{HM}(\Orth(L))=\frac{2}{g_{sp}^+(L)}\cdot |\det L|^{(\Lrank+1)/2}
\prod_{k=1}^{\Lrank}\pi^{-k/2}\Gamma(k/2)
 \prod_{p}\alpha_p(L)^{-1}
\end{equation}
where the $\alpha_p(L)$ are the local densities of the lattice $L$ and
$g_{sp}^+(L)$ is
the number of spinor genera in the genus of $L$.
\end{theorem}
\begin{proof}
This follows immediately from formulae
(\ref{eqn_Siegelvolume}), (\ref{eqn_gamma}), (\ref{eqn_HMvolume})
and (\ref{eqn_localglobal}).
\end{proof}

 \section{Applications}\label{Applications}

In this section we want to apply the above results to compute the
asymptotic behaviour of the dimension of spaces of cusp forms for a
number of specific groups. The main applications have to do with
locally symmetric varieties. In \cite{GHS1} we prove general type results
for the moduli spaces $\cF_{2d}$ of $\Kthree$ surfaces of
degree $2d$, but in that special case we can use a different
method. The results we have here are used in \cite{GHS2} to prove similar results
in greater generality.

\subsection{Groups}\label {subsec_groups}

We first have to clarify the various groups which will play a role. In
this section, $L$ will be an even indefinite lattice of signature
$(2,n)$, containing at least one hyperbolic plane as a direct summand.
By a classical result of Kneser we know that if the genus of an
indefinite lattice $L$ contains more than one class, then there is a
prime $p$ such that the quadratic form of $L$ can be diagonalised over the
$p$-adic numbers and the diagonal entries all involve distinct powers
of $p$ (see \cite[Chapter 15]{CS}). Therefore the genus of any
indefinite lattice with one hyperbolic plane contains only one class.

The group $\Orth(L)$ interchanges the two connected components of the domain
$\Omega_{L}$
and we define $\Orth^+(L)$ as the index $2$ subgroup which fixes each
of these components (as sets). For the connection with the spinor norm see
below.
The Hirzebruch--Mumford  volume of $\Orth^+(L)$ is twice that of $\Orth(L)$.
As an immediate corollary of Theorem  \ref{prop_mainformula} we obtain
\begin{theorem}\label{cor_mainformula}
Let $L$ be a lattice of signature $(2,n)$ ($n\ge 1$) containing at least
one hyperbolic plane.
Let $\Gamma$ be an arithmetic subgroup of $\Orth^+(L)$. Then
\begin{equation}\label{eqn_corollary}
\vol_{HM}(\Gamma)=2 \cdot [\POrth(L):\PGamma] |\det L|^{(n+3)/2}
\prod_{k=1}^{n+2}\pi^{-k/2}\Gamma(k/2)
 \prod_{p}\alpha_p(L)^{-1}.
\end{equation}
\end{theorem}

{\bf Remark.} In many interesting cases a subgroup $\Gamma$ is given
in terms of the orthogonal group of some sublattice $L_1$ of $L$.
In this case one can use the volume in order to calculate the index
(see Section~\ref{Tm} below).

Let $L^{\vee}= \Hom(L,\ZZ)$ be the dual lattice and
$A_L=L^{\vee}/L$. The finite group $A_L$ carries a discriminant
quadratic form $q_L$ with values in $\QQ/2\ZZ$ \cite [1.3]{N-int}. By
$\Orth(q_L)$ we denote the corresponding group of isometries and the
group $\widetilde{\Orth}(L)$, called the stable orthogonal group,
 is defined as the kernel of the natural
homomorphism $\Orth(L)\to \Orth(q_L)$. Since $L$ contains a hyperbolic
plane, it follows from \cite[Theorem 1.14.2]{N-int} that this map is
surjective.

The $(-1)$-{\em spinor norm}\/ on the group $\Orth(L \otimes \RR)$ can
be defined as follows. Every element $g$ can be represented as a
product of reflections
$$
g=\sigma_{v_1}\cdot \dots \cdot \sigma_{v_m}
$$
and, following Brieskorn \cite{Br}, we define
$$
\sn_{-1}(g)=
\left\{ \begin{array}{ll}
+1 &\mbox{ if } (v_k,v_k)>0\ \text{ for an even number of } v_k\\
-1 &\text { otherwise.}
\end{array}\right.
$$
This is independent of the representation of $g$ as a product of
reflections. We have already introduced the group $\Orth^+(L)$ as the
subgroup of $\Orth(L)$ leaving the two connected components of
$\Omega_L$ fixed. It is well known that
$$
\Orth^+(L)= \Ker(\sn_{-1})\cap \Orth(L).
$$
To see this, note that any reflection with respect to a vector of
negative square has $(-1)$-spinor norm equal to 1, and any reflection
with respect to a vector of positive square has $(-1)$-spinor norm equal
to $(-1)$ and interchanges the two components. Set
$$
\widetilde{\Orth}^+(L)= \widetilde{\Orth}(L) \cap \Orth^+(L).
$$
Finally the groups $\SO^+(L)$ and $\widetilde{\SO}^+(L)$ are defined
as the corresponding groups of isometries of determinant $1$.

\begin{lemma}
Let $N= |\Orth(q_L)|$. Then we have the following diagram of groups
with indices as indicated:
$$
\begin{matrix}
\widetilde{\Orth}(L)&\overset{N:1}{ \subset }& \Orth(L)\\
\bigcup\, {\scriptstyle{2:1}}&&\bigcup\, {\scriptstyle{2:1}}\\
\widetilde{\Orth}^+(L) &\overset{N:1}{ \subset }& \Orth^+(L)\\
\bigcup\, {\scriptstyle{2:1}}&&\bigcup\, {\scriptstyle{2:1}}\\
 \widetilde{\SO}^+(L) &\overset{N:1}{ \subset }& \SO^+(L).
\end{matrix}
$$
\end{lemma}
\begin{proof}
We shall first prove that the indices of the vertical inclusions are
all $2$. To do this, we choose a hyperbolic plane $U$ in $L$, which
exists by assumption. Let $e_1, e_2$ be a basis of $U$ with
$e_1^2=e_2^2=0$ and $e_1.e_2=1$. If $u=e_1-e_2$, $v=e_1+e_2$, then
$u^2=-2$, $v^2=2$ and the two reflections $\sigma_u$ and $\sigma_v$
belong to $\widetilde{\Orth}(L)$, since they act trivially on the
orthogonal complement of $U$. Moreover $\sn_{-1}(\sigma_v)=-1$ and
$\sn_{-1}(\sigma_u)=1$. Hence we can use $\sigma_v$ to conclude that
the top two vertical inclusions are of index $2$, whereas $\sigma_u$
shows the same for the bottom two vertical inclusions.

We have already observed that the natural map $\Orth(L)\to \Orth(q_L)$
is surjective, which shows that the top horizontal inclusion has index
$N$. Taking into account that the reflections $\sigma_u$ and
$\sigma_v$ act trivially on the discriminant form, we obtain that
$$
N=[\Orth(L):\widetilde{\Orth}(L)]=[\Orth^+(L):\widetilde{\Orth}^+(L)]=
[\SO^+(L):\widetilde{\SO}^+(L)].
$$
\end{proof}
Finally, we want to consider the projective groups $\POrth(L)$,
$\POrth^+(L)$ and $\POrthtd^+(L)$, i.e., the corresponding groups
modulo their centres. It follows immediately from the above diagram
that
\begin{equation}\label{eqn_index}
[\POrth(L):\POrthtd^+(L)]=
\left\{ \begin{array}{ll}
\ N &\mbox{ if } -\id \not \in \widetilde{\Orth}^+(L) \\
2N&\mbox{ if } -\id \in \widetilde{\Orth}^+(L).
\end{array}\right.
\end{equation}
Note that $-\id\in \widetilde{\Orth}^+(L)$ if and only if
$A_L$ is a $2$-group.

\subsection{Local densities}\label {subsec_densities}

Siegel's definition of local densities of a quadratic form over a
number field $K$ given by a matrix $S\in\Mat_{n\times n}(K)$ is
$$
\alpha_p(S)=\frac{1}{2}\lim_{r\to \infty}\
p^{-\frac{rn(n-1)}2}
|\left
\{X\in \Mat_{n\times n}(\ZZ_p)\hspace{-9pt}
\mod p^r;\, {}^tXSX\equiv S\hspace{-9pt}\mod p^r\right\}|.
$$

The local densities can be calculated explicitly, at least
in the cases where $K=\QQ$ or a quadratic extension of $\QQ$
(see chapter~5 of the book \cite{Ki} and references there).
For the convenience of the reader we include the formulae over
$\QQ$ in the present paper. To calculate $\alpha_p(L)$ one should know
the Jordan decomposition of the lattice $L$ over the local ring
$\ZZ_p$ of $p$-adic integers. The main difficulties arise for $p=2$:
see \cite[Theorem 5.6.3]{Ki}.

Let us introduce some notation. Let $L$ be a $\ZZ_p$-lattice in a
regular (i.e. nondegenerate) quadratic space over $\QQ_p$
of rank $n$, and let $(\bdv_i)$ be a basis of $L$. There are
two invariants of $L$:
the scale
$$
\scale(L)=\{(\bdx,\bdy)_L\,;\,\bdx,\bdy\in L\}
$$
and the norm
$$
\norm(L)=\{{\textstyle\sum} a_\bdx(\bdx,\bdx)_L\,;\, \bdx\in L,\
a_\bdx\in \ZZ_p\, \}.
$$
We have $2\scale(L)\subset \norm(L) \subset \scale(L)$.
In fact, over $\ZZ_p$ ($p\ne 2$) we have $\norm(L)=\scale(L)$,
whereas over $\ZZ_2$ we have either $\norm(L)=\scale(L)$ or
$\norm(L)=2\scale(L)$.

$L$ is called $p^r$-modular, for $r\in\ZZ$, if the matrix $ p^{-r}(\bdv_i,\bdv_j)_L$
belongs to $\GL_n(\ZZ_p)$. In this case we can write $L$ as
the scaling $N(p^r)$ of a unimodular lattice $N$. By a hyperbolic
space we mean a (possibly empty) orthogonal sum of hyperbolic planes.

A regular lattice $L$ decomposes as the orthogonal sum of lattices
$\bigoplus_{j\in\ZZ}L_j$, where $L_j$ is a $p^j$-modular lattice of rank
$n_j\in \ZZ_{\ge 0}$. Put
$$
w=\sum_j  jn_j\bigl( (n_j+1)/2 + \sum_{k>j} n_k \bigr)
$$
and
$$
P_p(n)=\prod_{i=1}^n (1-p^{-2i}).
$$
For a regular quadratic space $W$ over the finite field $\ZZ/p\ZZ$ one
puts
$$
\chi(W)=\begin{cases}
\hphantom{-}0&\quad \text{if $\dim W$ is odd,}\\
\hphantom{-}1&\quad \text{if $W$ is a hyperbolic space,}\\
-1&\quad \text{otherwise.}
\end{cases}
$$
For a unimodular lattice $N$ over $\ZZ_2$ with $\norm(N)=2\scale(N)$
we define $\chi(N)=\chi(N/2N)$, where $N/2N$ is given the structure of
a regular quadratic space over $\ZZ/2\ZZ$ via the quadratic form
$Q(\bdx)=\frac{1}{2}(\bdx,\bdx)_N\bmod 2$.
\medskip

For the local density $\alpha_p(L)$ for $p\ne 2$ we have the formula
\begin{equation}\label{eqn_ap}
\alpha_p(L)=2^{s-1}p^wP_p(L)E_p(L)
\end{equation}
where $s$ is the number of non-zero $p^j$-modular terms $L_j$
in the orthogonal decomposition of $L$, and
$$
P_p(L)=\prod_j P_p([n_j/2]), \qquad
E_p(L)=\prod_{j,\, L_j\ne 0} \bigl(1+\chi(N_j)p^{-n_j/2}\bigr)^{-1}
$$
where $L_j$ is the $p^j$-scaling of the unimodular lattice $N_j$ and $[n_j/2]$ denotes the integer part.
\medskip

The local density $\alpha_2(L)$ is given by
\begin{equation}\label{eqn_a2}
\alpha_2(L)=2^{n-1+w-q}P_2(L)E_2(L).
\end{equation}
In this formula $q=\sum_j q_j$ where
$$
q_j=
\begin{cases}
0&\quad \text{if $N_j$ is even,}\\
n_j&\quad \text{if $N_j$ is odd and $N_{j+1}$ is even,}\\
n_j+1&\quad \text{if $N_j$ and $N_{j+1}$ are odd}.
\end{cases}
$$
A unimodular lattice $N$ over $\ZZ_2$ is even if it is trivial or
if $\norm(N)=2\ZZ_2$, and odd otherwise. Any unimodular lattice can be
represented as the orthogonal sum $N=N^{\even}\oplus N^{\odd}$ of even and odd
sublattices such that $\rank N^{\odd}\le 2$. Then we put
$$
P_2(L)=\prod_j P_2(\rank N^{\even}_j/2).
$$
The second factor is $E_2(L)=\prod_j E_j^{-1}$, where $E_j$ is defined by
$$
E_j=\frac{1}{2}(1+\chi(N_j^{\even})2^{-\rank N_j^{\even}/2})
$$
if both $N_{j-1}$ and $N_{j+1}$ are even, unless
$N_j^{\odd}\cong \latt{\epsilon_1}\oplus \latt{\epsilon_2}$ with
$\epsilon_1\equiv\epsilon_2\mod 4 $: in all other cases we put
$E_j=1/2$.

We note that $E_j$ depends on $N_{j-1}$, $N_{j}$ and $N_{j+1}$
and $E_j=1$ if all of them are trivial.
Also $q_j$ depends on $N_j$ and $N_{j+1}$ and
$q_j=0$ if $N_j$ is trivial.

\subsection{The even unimodular lattices \label{evenunimodular}
 ${II}_{2,8m+2}$}\label{sub_evenunimodular}

We start with the example
$$
{II}_{2,8m+2}=2U \oplus m E_8(-1), \qquad \text{where}\quad m\ge 0
$$
which is a natural series of even unimodular lattices of signature
$(2,8m+2)$. Note that ${II}_{2,26} \cong 2U \oplus \Lambda$, where $\Lambda$
is the Leech lattice.

The local densities are easy to calculate, since for every prime $p$
the lattice ${II}_{2,8m+2} \otimes \ZZ_p$ over the $p$-adic integers is a
direct sum of hyperbolic planes. Then using (\ref{eqn_ap}) and
(\ref{eqn_a2}) we obtain
$$
\alpha_{p}({II}_{2,8m+2})=2^{\delta_{2,p}(8m+4)}P_p(4m+2)(1+p^{-(4m+2)})^{-1}
$$
where $\delta_{2,p}$ is the Kronecker delta. By our main formula
(\ref{eqn_mainvolume}) we obtain
$$
\vol_{HM}({\Orth}^+({II}_{2,8m+2}))=
2^{-(8m+2)}\gamma_{8m+4}
\zeta(2)\zeta(4)\cdot\ldots\cdot \zeta (8m+2)\zeta(4m+2)
$$
where $\gamma_{8m+4}$ is as in formula (\ref{eqn_gamma}). In order
to simplify this expression we use the $\zeta$-identity
\begin{equation}\label{E:val-zeta}
\pi^{-\frac{1}2-2k}\Gamma(k)\Gamma\left(k+\frac{1}2\right)\zeta(2k)=
(-1)^k\zeta(1-2k)=(-1)^{k+1}\frac{B_{2k}}{2k}.
\end{equation}
Together with
\begin{eqnarray*}
\lefteqn{\pi^{-(4m+2)}\Gamma(4m+2)\zeta(4m+2)}\hspace{15pt}\\
&=& 2^{4m+1}\pi^{-\frac{1}{2}-(4m+2)}
\Gamma(2m+1)\Gamma(\frac{4m+3}{2})\zeta(4m+2)\\
&=& 2^{4m+1}\frac{B_{4m+2}}{4m+2}
\end{eqnarray*}
where the first equality comes from the Legendre duplication formula
of the $\Gamma$-function, and the second equality is again a
consequence of the $\zeta$-identity, we obtain
$$
\vol_{HM}(\Orth^+({II}_{2,8m+2}))= 2^{- (4m+1)}
\frac{B_2\cdot B_4\cdot\ldots\cdot B_{8m+2}}{(8m+2)!!}
\cdot \frac{B_{4m+2}}{4m+2}.
$$
Here $(2n)!!=2 \cdot 4 \ldots 2n$. Since the discriminant group of
the lattice ${II}_{2,8m+2}$ is trivial, we have the equality
$$
\vol_{HM}(\widetilde{\Orth}^+ ({II}_{2,8m+2}))=\vol_{HM}(\Orth^+({II}_{2,8m+2})).
$$
 In a similar way one can derive a formula for any indefinite
unimodular lattice of signature $(r,s)$. For example, for the odd unimodular
lattice $M$ defined by $x_1^2+\dots+x_r^2-x_{r+1}^2-\dots - x_{r+s}^2$
we have to take into account that the even $(M\otimes \ZZ_2)^{\even}$ and odd
$(M\otimes \ZZ_2)^{\odd}$ parts  of the lattice
$M$ over $2$-adic numbers depend on $r+s\mod 2$ and $r-s\mod 8$
(see \cite{BG} for a different approach  in this special case).

We can now use this to compute dimensions of cusp forms for this
group and we obtain
\begin{eqnarray*}
\lefteqn{\dim S_k(\widetilde{\Orth}^+ ({II}_{2,8m+2}),\det{}^{\varepsilon})=}\\ 
&&\frac{2^{-4m}}{(8m+2)!}\cdot
\frac{B_2\cdot B_4\cdot\ldots\cdot B_{8m+2}}{(8m+2)!!}
\cdot \frac{B_{4m+2}}{4m+2}k^{8m+2}+O(k^{8m+1}).
\end{eqnarray*}
Here $\varepsilon= \pm 1$ and we must assume that $k$ is even, since otherwise there are no forms
for trivial reasons.

\subsection{The lattices $T_{2,8m+2}$}\label{Tm}

The orthogonal group of the lattice ${II}_{2,8m+2}$ for $m=2$ defines
an irreducible component of the branch divisor of the modular variety
$\cF_{2d}^{(m)}$. The same branch divisor contains
another component defined by the lattice
$$
T_{2,8m+2}=U\oplus U(2)\oplus mE_8(-1)
$$
of discriminant $4$. We note that this lattice is not maximal.
 For a prime number $p\ne 2$ the $p$-local densities
of the lattice $T$ and $M$ coincide. Let us calculate $\alpha_2(T)$.
Over the $2$-adic ring we have $T_{2,8m+2}\otimes \ZZ_2\cong
 (4m+1)U\oplus U(2)$.
We have (see (\ref{eqn_a2}))
$$
N_0=N_0^{\even}=(4m+1)U,\quad N_1=N_1^{\even}=U,\quad w=3,\quad q=0,
$$
$$
E_0=\frac1{2}(1+2^{-(4m+1)}),\qquad \ E_1=\frac1{2}(1+2^{-1}).
$$
Thus
$$
\alpha_2(T_{2,8m+2})=2^{8m+7}(1-2^{-2})\cdot \dots \cdot
(1-2^{-8m})(1-2^{-(4m+1)}).
$$
We note that $[\POrth^+(T_{2,8m+2}): \POrthtd^+(T_{2,8m+2})]=2$
since the finite orthogonal discriminant group of $T_{2,8m+2}$
is isomorphic to $\ZZ/2\ZZ$. As a result we get
\begin{eqnarray*}
\lefteqn{\vol_{MH}(\widetilde{\Orth}^+(T_{2,8m+2}))=}\\
&&
2\,\gamma_{8m+4}\zeta(2)\cdot\dots\cdot \zeta(8m+2)\zeta(4m+2)
(1+2^{-(4m+1)})(1-2^{-(4m+2)}).
\end{eqnarray*}
Using the formula for the volume of ${II}_{2,8m+2}$ we see that
\begin{equation}\label{compareMT}
\frac{\vol_{MH}\widetilde{\Orth}^+(T_{2,8m+2})}
{\vol_{MH}\widetilde{\Orth}^+({II}_{2,8m+2})}=(2^{4m+1}+1)(2^{4m+2}-1).
\end{equation}

If $L_1$ is a sublattice of finite index of a lattice $L$
then $\widetilde\Orth^+(L)$ is a subgroup of $\widetilde\Orth^+(L_1)$.
One can use the formula of Theorem~\ref{prop_mainformula} to calculate
easily the index $[\widetilde\Orth^+(L_1):\widetilde\Orth^+(L)]$.
For example, formula~(\ref{compareMT}) above gives the index of
$\widetilde\Orth^+({II}_{2,8m+2})$ in $\widetilde\Orth^+(T_{2,8m+2})$.
This method is much shorter than the calculation in terms of finite
geometry over $\ZZ/2\ZZ$.

\subsection{The lattices $L_{2d}^{(m)}$}\label{sub_L2d}

We consider the lattice
$$
L_{2d}^{(m)}=2U\oplus mE_8(-1)\oplus \latt{-2d}
$$
of signature $(2,8m+3)$. The lattice $L_{2d}^{(m)}$ is not maximal if
$d$ is not square free.
This lattice is of particular interest, as
the lattice $L_{2d}^{(2)}$ is closely related to the moduli space of
polarised $\Kthree$ surfaces of degree $2d$. More precisely, the
quotient space
$$
{\cF}_{2d}= \widetilde{\Orth}^+(L_{2d}^{(2)}) \backslash \cD_{L_{2d}^{(2)}}
$$
is the moduli spaces of $\Kthree$ sufaces of degree $2d$. As we shall
see, there is also a relation to Siegel modular forms for both the
group $\Sp(2,\ZZ)$ and the paramodular group.

Again, the lattices over the $p$-adic integers are easy to understand,
since $E_8(-1) \otimes \ZZ_p$ is the direct sum of four copies of a
hyperbolic plane. By (\ref{eqn_ap}) and (\ref{eqn_a2}) we find
$$
\begin{aligned}
\alpha_{p}(L_{2d}^{(m)})&=P_p(4m+2)
\quad&&\text{ if } p\not| 2d\\
\alpha_{p}(L_{2d}^{(m)})&=2p^sP_p(4m+2)(1+p^{-(4m+2)})^{-1}
\quad&&\text{ if $p$ is odd, } p^s\| d\\
\alpha_{2}(L_{2d}^{(m)})&=2^{8m+6}P_2(4m+2)
\quad&&\text{ if $d$ is odd}\\
\alpha_{2}(L_{2d}^{(m)})&=2^{8m+7+s}P_2(4m+2)(1+2^{-(4m+2)})^{-1}
\quad&&\text{ if $d$ is even, } \ 2^s\|d.
\end{aligned}
$$
Therefore
$$
\prod_{p}\alpha_p(L^{(m)}_{2d})^{-1}=
\zeta(2)\zeta(4) \dots \zeta(8m+4)\,
(2d)^{-1}2^{-\rho(d)-8m-5}\prod_{p|d}(1+p^{-(4m+2)})
$$
where $\rho(d)$ denotes the number of prime divisors of $d$.

We shall need the following.
\begin{lemma}
\label{lem:discriminant0}
Let
$R= \latt{-2d}$.
Then the order of the discriminant group $\Orth(q_R)$ is $2^{\rho(d)}$.
\end{lemma}
\begin{proof}
Let $g$ be the standard generator of $A_R=\ZZ / 2d\ZZ$, given by the
equivalence class of $1$. Then $q_R(g)= -1/2d \mod 2\ZZ$. If $\varphi
\in \Orth(q_R)$, then $\varphi(g)=xg$ for some $x$ with
$(x,2d)=1$. Hence $\varphi$ is orthogonal if and only if
$$
- \frac{x^2}{2d} \equiv -\frac{1}{2d} \mod 2\ZZ,
$$
or equivalently
$$
x^2 \equiv 1 \mod 4d\ZZ.
$$ It is not difficult to check that this equation has $2^{\rho(d)+1}$
solutions modulo $4d\ZZ$, and hence $2^{\rho(d)}$ solutions modulo
$2d\ZZ$.
\end{proof}
{}From this it follows also that $|A_{L^{(m)}_{2d}}|=2^{\rho(d)}$.

{}From (\ref{eqn_index}) it follows that
$$
[\POrth(L_{2d}^{(m)}):\POrthtd^+(L_{2d}^{(m)})]=2^{\rho(d)} \quad
\mbox{ if } \quad d>1
$$
and $2$ if $d=1$. We first assume that $d>1$. We put $n=8m+3$,
which is the dimension of the homogeneous domain. It follows from
Corollary~\ref{cor_mainformula} that
$$
\vol_{HM}(\widetilde{\Orth}^+(L_{2d}^{(m)}))=2^{\rho(d)+1}(2d)^{\frac{n+3}{2}}
\gamma_{n+2}\prod_{p}\alpha_p(L^{(m)}_{2d})^{-1}.
$$
If $d=1$ we have to multiply this formula by a factor $2$. Using
the $\zeta$ identity, a straightforward calculation gives (again for
$d>1$ and $n=8m+3$)
$$
\vol_{HM}(\widetilde{\Orth}^+(L_{2d}^{(m)}))=
\left( \frac{d}{2} \right )^\frac{n+1}{2}
\prod_{p|d}(1+p^{-\frac{n+1}{2}})\cdot
\frac{|B_2\cdot B_4\cdot\dots\cdot B_{n+1}|}
{(n+1)!!}.
$$

We want to apply this to the moduli space of $\Kthree$ surfaces
of degree $2d$. This is the case $m=2$: the dimension of the domain is
$n=19$. Using Hirzebruch-Mumford proportionality and specialising the
above volume computation to this case, we compute the dimension of the
spaces of cusp forms:
\begin{eqnarray*}
\lefteqn{\dim S_{k}(\widetilde{\Orth}^+(L_{2d}^{(2)}), \det{}^{\varepsilon}) =} \\
&& \frac{2^{-9}}{19!}d^{10} \cdot \prod_{p|d}(1+p^{-10})
\frac{|\,B_2\cdot B_4\cdot\ldots\cdot B_{20}\,|}{20!!}\cdot k^{19}+O(k^{18})
\end{eqnarray*}
which holds for $d>1$, with an additional factor $2$ for
$d=1$. In the latter case we must assume that $k$ and $\varepsilon$ have the same parity.
For $d>1$ there is no restriction since $- \id \notin \widetilde{\Orth}^+(L_{2d}^{(2)})$.
This should be compared to 
Kondo's formula \cite{Ko} where, however, the Hirzebruch-Mumford
volume has not been computed explicitly. It should also be noted that
Kondo uses the geometric, rather than the arithmetic, weight.

\subsubsection{Siegel modular forms}

The case $m=0$ gives applications to Siegel modular forms. We shall
first consider the case $d=1$. Recall that
$$
\widetilde{\SO}^+(L_2^{(0)}) \cong \widetilde{\Orth}^+(L_2^{(0)})/\{\pm \id\} \cong
\Sp(2,\ZZ)/\{\pm \id\}.
$$
{}From our previous computation we obtain that
$$
\vol_{HM}(\widetilde{\SO}^+(L_{2}^{(0)}))=
\vol_{HM}(\widetilde{\Orth}^+(L_{2}^{(0)}))=2^{-4}|B_2B_4|
$$
and by Hirzebruch-Mumford proportionality this gives
$$
\dim S_k(\Sp(2,\ZZ))=2^{-4}3^{-1}|B_2B_4|k^3 + O(k^2).
$$
Note that this coincides with \cite[p. 428]{T}, taking into account
that Tai's formula refers to modular forms of weight $3k$. Tai uses
Siegel's computation of the volume of the group $\Sp(2,\ZZ)$, rather
than the orthogonal group.

\subsubsection{The paramodular group}

Finally, we consider the case $m=0$ and $d > 1$. This is closely
related to the so-called {\em paramodular\/} group $\Gamma_{d}^{(\Sp)}$,
which gives rise to the moduli space of $(1,d)$-polarised abelian
surfaces. In fact
$$
\widetilde{\SO}^+(L_{2d}^{(0)}) \cong \PGamma_{d}^{(\Sp)}
$$
by \cite[Proposition 1.2]{GH}. We note that in this case
$$
[\widetilde{\Orth}^+(L_{2d}^{(0)}):\widetilde{\SO}^+(L_{2d}^{(0)})]=2
$$
and that $-\id$ is in neither of these groups. Hence
\begin{eqnarray*}
\vol_{HM}(\widetilde{\SO}^+(L_{2d}^{(0)})) &=&
2\vol_{HM}(\widetilde{\Orth}^+(L_{2d}^{(0)}))\\
&=&
2^{-4}d^2\prod_{p|d}(1+p^{-2})|B_2B_4|
\end{eqnarray*}
and by Hirzebruch-Mumford proportionality
$$
\dim S_{k}(\Gamma_{d}^{(\Sp)})=\frac{d^2}{3\cdot
 2^4}\prod_{p|d}(1+p^{-2})|B_2B_4| k^3+O(k^2).
$$
This agrees with \cite[Proposition 2.2]{Sa}, where this formula was
derived for $d$ a prime.

\subsection{Lattices associated to Heegner divisors} \label{Heegner}
\label{subsec_rank20lattices}

We shall conclude this section by computing the volume of two lattices
of rank $8m+4$. Both of these lattices  $K_{2d}^{(m)}$ and $N_{2d}^{(m)}$
arise from the  $(-2)$-reflective part of the ramification
divisor of the quotient map
$$
\cD_{L_{2d}^{(m)}} \to \widetilde{\Orth}^+_{L_{2d}^{(m)}} \backslash
\cD_{L_{2d}^{(m)}}=\cF_{2d}^{(m)}.
$$
For $m=2$ this is the moduli space of $\Kthree$ surfaces of degree
$2d$. For $m=0$ and a prime $d$ we get the moduli of Kummer surfaces
associated to $(1,d)$-polarised abelian surfaces
(see \cite{GH}).
Since the branch locus of the quotient map gives rise to obstructions for extending
pluricanonical forms defined by modular forms, knowledge of their volumes is important
for the computation of the Kodaira dimension of $\cF_{2d}^{(m)}$.

\subsubsection{The lattices $K_{2d}^{(m)}$} \label{subsubsec_volumeK}

We consider the lattice
$$
K_{2d}^{(m)} = U \oplus mE_8(-1) \oplus \latt{2} \oplus \latt{-2d}
$$
where $d$ is a positive integer. We first have to determine the
local densities for this lattice. Since $\det(K_{2d}^{(m)})=4d$, this
lattice is equivalent to the following lattices over the $p$-adic
integers for odd primes $p$:
$$
K_{2d}^{(m)} \otimes \ZZ_p \cong (4m+1)\,U\oplus
\begin{cases}
U\quad&\text{if } \left(\frac{4d}{p}\right)=1\\
x^2-4dy^2\quad&\text{if } \left(\frac{4d}{p}\right)=-1.
\end{cases}
$$
For the local densities we obtain from equations (\ref{eqn_ap}) and
(\ref{eqn_a2})
$$
\begin{aligned}
\alpha_{p}(K_{2d}^{(m)})&=P_p(4m+1)(1-\left(\frac{4d}{p}\right)p^{-(4m+2)})
\quad&&\text{ if }\ p\not| d\\
\alpha_{p}(K_{2d}^{(m)})&=2p^sP_p(4m+1)
\quad&&\text{ if }\ p^s\|d\\
\alpha_{2}(K_{2d}^{(m)})&=2^{8m+v(d)}P_2(4m+1)
\end{aligned}
$$
where
$v(d)=6$ if $d\equiv 1 \mod 4$,
$v(d)=7$ if $d\equiv -1 \mod 4$,
$v(d)=8$ if $d\equiv 2 \mod 4$, and
$v(d)=8+s$ if $d\equiv 0 \mod 4$ and $2^s\,||\,d$.

{}From this we obtain that
\begin{equation} \label{eqn_prodalphaK}
\prod_{p}\alpha_p(K_{2d}^{(m)})^{-1}= A_2(d) d^{-1}
\zeta(2)\zeta(4)\dots\zeta(8m+2)
L(4m+2,\left(\frac{4d}{*}\right)),
\end{equation}
where
$$
A_2(d)=\begin{cases}
2^{-\rho(d)-8m-6}\quad&\text{if } d\equiv 1,\,2 \mod 4\\
2^{-\rho(d)-8m-7}\quad&\text{if } d\equiv 0,\,3 \mod 4.
\end{cases}
$$
Application of our main formula (\ref{eqn_mainvolume}) then gives
\begin{equation} \label{eqn_vol1K}
\vol_{HM}(\Orth^+(K_{2d}^{(m)}))=
4 \cdot (4d)^{\frac{8m+5}{2}} \cdot \prod_{k=1}^{8m+4}
\pi^{-\frac{k}{2}}\Gamma(\frac{k}{2}) \cdot \prod_{p}\alpha_p(K_{2d})^{-1}.
\end{equation}
Combining formulae (\ref{eqn_prodalphaK}) and (\ref{eqn_vol1K}) and the
$\zeta$-identity (\ref{E:val-zeta}) leads to
\begin{eqnarray}
\lefteqn{\vol_{HM}(\Orth^+(K_{2d}^{(m)}))=}\nonumber \\
&&
C_2(d) d^{{\frac{8m+3}{2}}}\pi^{-(4m+2)}
\Gamma(4m+2)L(4m+2,\left(\frac{4d}{*}\right))
\frac{\,B_2B_4\ldots B_{8m+2}\,}{(8m+2)!!}\nonumber
\end{eqnarray}
where
$$
C_2(d)=\begin{cases}
2^{-\rho(d)+1}\quad&\text{if } d\equiv 1,\,2 \mod 4\\
2^{-\rho(d)}\quad&\text{if } d\equiv 0,\,3 \mod 4.
\end{cases}
$$
For applications it is also important to compute the volume with
respect to the group $\widetilde{\Orth}^+(K_{2d}^{(m)})$.
For this, we have to know the order of the group of isometries
of the discriminant group.
\begin{lemma}
\label{lem:discriminant}
Let $S=\latt{2}\oplus \latt{-2d}$. The order of the discriminant group is
$$
|\Orth(q_S)|=
\left\{ \begin{array}{ll}
2^{1+\rho(d)} & \mbox{ if } d\equiv -1 \mod 4 \mbox{ or } d \mbox{ is
 divisible by } 8\\
2^{\rho(d)} & \mbox{ for all other } d.
\end{array}\right.
$$
\end{lemma}
\begin{proof}
We denote the standard generators of $\ZZ / 2d\ZZ$ and $\ZZ / 2\ZZ$ by
$g$ and $h$ respectively. We shall first consider automorphisms
$\varphi$ with $\varphi(g)=xg$. Then orthogonality implies $x^2 \equiv
1 \mod 4d\ZZ$ which means, in particular, that $x$ is odd and
$(x,2d)=1$. We then have $\varphi(dg)=dg$. We cannot have that
$\varphi(h)=dg+h$, because orthogonality implies that for the bilinear
form $B_q$, defined by the quadratic form $q=q_S$, we have
$B_q(xg,dg+h)=B_q(g,h)=0$ and hence $-x/2 \equiv 0 \mod \ZZ$, which
shows that $x$ is even, a contradiction. Hence $\varphi(h)=h$ and
$\varphi = \varphi' \times \operatorname{id}$ where $\varphi' \in
\Orth(q_R)$ (with $R= \latt{-2d}$). In this way we obtain
$2^{\rho(d)}$ elements in $\Orth(q_S)$.

We shall now investigate automorphisms with $\varphi(g)=xg+h$.
Then $q(g)=q(\varphi(g))$ implies the condition
$$
x^2 \equiv 1+d \mod 4d\ZZ.
$$
It is not hard to check that this only has solutions if either $d
\equiv -1 \mod 4$ or $d$ is divisible by $8$. We shall distinguish
between the cases $d$ even and $d$ odd. In the first case $x$ must be
odd and $(x,2d)=1$. Moreover $\varphi(dg)=dg$ and the only possibility
for an orthogonal automorphism is $\varphi(h)=dg+h$ and indeed this
gives rise to another $2^{\rho(d)}$ orthogonal automorphisms. Now
assume $d$ is odd. Then $x$ is even and $(x,d)=1$. In this case
$\varphi(dg)=h$ and the only possibility to obtain an orthogonal
automorphism is $\varphi(h)=dg$. Once more, this gives another
$2^{\rho(d)}$ orthogonal automorphims and this proves the lemma.
\end{proof}
By formula (\ref{eqn_index}) it then follows that
$$
[\POrth(K_{2d}^{(m)}):\POrthtd^+(K_{2d}^{(m)})]=
\begin{cases}
2\quad&\text{if }d=1\\
2^{\rho(d)}
\quad&\text{if }d\equiv 1,\,2 \mod 4,\ d>1\\
2^{\rho(d)+1}
\quad&\text{if }d\equiv 3 \mod 4.
\end{cases}
$$
Therefore
\begin{eqnarray}
\lefteqn{\vol_{HM}(\Tilde{\Orth}^+(K_{2d}^{(m)}))=
2^{\delta_{1,d}-\delta_{4,d(8)}}\frac{\,B_2B_4\ldots B_{8m+2}\,}{(8m+2)!!}}
\nonumber\\
\label{eqn_vol2K}&&
\cdot d^{\frac{8m+3}{2}}
\pi^{-(4m+2)}
\Gamma(4m+2)L(4m+2,\left(\frac{4d}{*}\right))
\end{eqnarray}
where $d(8)$ denotes $d\mod 8$ and $\delta_{*,*}$  is the Kronecker symbol.

We want to reformulate this result in terms of generalised
Bernoulli numbers. In order to avoid too many different cases, we
restrict here to $ d \not\equiv 0 \mod 4$ (but it is clear how to
remove this restriction). If $d=d_0t^2$, with $d_0$ a positive and
square-free integer, then the discriminant of the real quadratic field
$\QQ(\sqrt{d})$ is equal to
$$
D=\begin{cases}
d_0\quad&\text{if } d\equiv 1 \mod 4\\
4d_0\quad&\text{if } d\equiv 2,\,3\mod 4.
\end{cases}
$$
Note that
\begin{equation}\label{eqn_discr}
d^{\frac{8m+3}{2}}=t^{8m+3}D^{\frac{8m+3}{2}} \cdot\begin{cases}
1 \quad&\text{if } d\equiv 1\ \ \mod 4\\
2^{-(8m+3)} \quad&\text{if } d\equiv 2,\,3\mod 4.
\end{cases}
\end{equation}
Let $\chi_D$ be the quadratic character of this field.
Then
\begin{equation}\label{eqn_Lfunction}
L(s, \left(\frac{4d}{*}\right))=L(s,\chi_D)
\prod_{p|2t}(1-\chi_D(p)p^{-s}).
\end{equation}
The character $\chi_D$ is an even primitive character modulo $D$, and
the Dirichlet $L$-function $L(s,\chi_D)$ satisfies the functional
equation
\begin{equation}\label{eqn_DirichletLfunctional}
\pi^{-\frac{s}2}\Gamma(\frac{s}2)D^sL(s,\chi_D)
=\pi^{-\frac{1-s}2}\Gamma(\frac{1-s}2)D^{\frac{1}2}L(1-s,\chi_D).
\end{equation}
Moreover
$$
L(1-k,\chi_D)=-\frac{B_{k,\chi_D}}{k}
$$
where $B_{k,\chi_D}$ is the corresponding generalised Bernoulli number.
Using the functional equation (\ref{eqn_DirichletLfunctional}) we obtain
\begin{eqnarray}\label{eqn_Gammatransform}
\pi^{-(4m+2)}\Gamma(4m+2)D^{\frac{8m+3}2}L(4m+2,\chi_D)&=&
-2^{4m+1}L(1-(4m+2),\chi_D)\nonumber\\
&=&2^{4m+1}\frac{B_{4m+2,\chi_D}}{4m+2}.
\end{eqnarray}
Combining (\ref{eqn_vol2K}), (\ref{eqn_discr}),
(\ref{eqn_Lfunction}), (\ref{eqn_Gammatransform}) and the result of
Lemma \ref{lem:discriminant} then gives the
result
\begin{eqnarray}
\lefteqn{\vol_{HM}(\Tilde{\Orth}^+(K_{2d}^{(m)}))=}\nonumber\\
\label{eqn_vol3K}&&
F_2(d) t^{8m+3}\frac{\,B_2B_4\ldots B_{8m+2}\,}{(8m+2)!!}
\frac{B_{4m+2,\chi_D}}{4m+2}
\prod_{p|2t}(1-\chi_D(p)p^{-(4m+2)})
\end{eqnarray}
where
$$
F_2(d)=\begin{cases}
2^{4m+2} \quad&\text{if } d\equiv 1\mod 4\\
2^{-4m-1} \quad&\text{if } d\equiv 2,\ 3\mod 4.
\end{cases}
$$
Using this, together with Hirzebruch-Mumford proportionality, we
finally find
that  $\dim S_k(\widetilde{\Orth}^+(K_{2d}^{(m)}))$ grows as
$$
\frac{G_2(d)}{(8m+2)!}
\frac{\,B_2\cdot B_4 \dots B_{8m+2}\,}{(8m+2)!!}\cdot
\frac{B_{4m+2,\chi_D}}{4m+2}\,
t^{8m+3}\prod_{p|2t}(1-\chi_D(p)p^{-(4m+2)})\,k^{8m+2}
$$
where
$$
G_2(d)=
\begin{cases}
2^{4m+2+\delta_{1,d}}\quad&\text{if }d\equiv 1\mod 4,\\
2^{-(4m+1)}\quad&\text{if }d\equiv 2,\ 3\mod 4.
\end{cases}
$$

\subsubsection{The lattices $N_{2d}^{(m)}$} \label{subsubsec_volumeN}

We assume that $d \equiv 1 \mod 4$ and consider the even lattice
$$
N_{2d}^{(m)}=U \oplus mE_8(-1) \oplus
\begin{pmatrix}
2&1\\1&\frac{1-d}{2}
\end{pmatrix}.
$$
We first have to understand this lattice over the $p$-adic
integers. If $p > 2$ then $2$ is a $p$-adic integer and we
have the following equality for the anisotropic binary form in
$N_{2d}^{(m)}$:
$$
\frac{1-d}{2}x^2 + 2xy + 2y^2= -\frac{d}{2} x^2 + 2(y+ \frac{x}{2})^2.
$$
Depending on whether $d$ is a square in $\ZZ_p^*$ or not,
we then obtain from the classification theory of quadratic forms
over $\ZZ_p$ that
$$
N_{2d}^{(m)}\otimes \ZZ_p \cong (4m+1)\, U\oplus
\begin{cases}
U\quad&\text{if } \left(\frac{d}{p}\right)=1\\
-dx^2 + y^2\quad&\text{if } \left(\frac{d}{p}\right)=-1.
\end{cases}
$$
We now turn to $p=2$. Recall that there  are only two even unimodular
binary forms over $\ZZ_2$, namely the hyperbolic plane and the form
given by the matrix $\begin{pmatrix} 2&1\\1& 2
\end{pmatrix}$. This implies that
$$
N_{2d}^{(m)}\otimes \ZZ_2 \cong (4m+1)\, U \oplus
\begin{cases}
\quad U\quad&\text{if } d\equiv 1\mod 8\\
\begin{pmatrix}
2&1\\1&2
\end{pmatrix}
\quad&\text{if } d\equiv 5\mod 8.
\end{cases}
$$
Once again by (\ref{eqn_ap}) and (\ref{eqn_a2}) we find for the local
densities that
$$
\begin{aligned}
\alpha_{p}(N_{2d}^{(m)})&=P_p(4m+1)(1-\left(\frac{d}{p}\right)p^{-(4m+2)})
\quad&&\text{ if } p\not| d\\
\alpha_{p}(N_{2d}^{(m)})&=2p^sP_p(4m+1)
\quad&&\text{ if }p^s\|d\\
\alpha_{2}(N_{2d}^{(m)})&=2^{8m+4}P_2(4m+1)(1-\left(\frac{d}{2}\right)2^{-(4m+2)}).
\end{aligned}
$$
We are interested mainly in the group
$\widetilde{\Orth}^+(N_{2d}^{(m)})$. For this we need the next lemma.
\begin{lemma}
\label{lem:discriminant2}
Let
$$
T=\begin{pmatrix}
2&1\\1&\frac{1-d}{2}
\end{pmatrix}
$$
Then $A_T \cong \ZZ / 2d\ZZ$ and
$$
|\Orth(q_T)|= 2^{\rho(d)}.
$$
\end{lemma}
\begin{proof}
Since $\det(T)=-d$, the discriminant group has order $d$. In fact, it
is cyclic of order $d$. To see this, let $e$ and $f$ be the basis
with repect to which the form is given by the matrix $T$. Then
$(e-2f)/d$ is in the dual lattice and its class, say $h$, generates
the group $A_T$. Every homomorphism of $A_T$ is of the form
$\varphi(h)=xh$, and it is an isometry if and only if $x^2 \equiv 1
\mod 2d$. This equation has $2^{\rho(d)}$ solutions modulo $d \ZZ$.
\end{proof}

It now follows from (\ref{eqn_index}) that
$$
[\POrth(N_{2d}^{(m)}):\POrthtd^+(N_{2d}^{(m)})]=
\begin{cases}
2^{\rho(d)}
\quad&\text{if }\ d\equiv 1 \mod 4 \text{ and } d\ne 1\\
2\quad&\text{if }\ d=1.
\end{cases}
$$
By the same calculation as in the preceding example we find now that
\begin{eqnarray}
\lefteqn{
\vol_{HM}(\Tilde{\Orth}^+(N_{2d}^{(m)}))=
2^{\delta_{1,d}-8m-3}
\frac{B_2 B_4 \dots B_{8m+2}\,}{(8m+2)!!}
\cdot}\nonumber\\
\label{eqn_vol2N}&&
d^{\frac{8m+3}2}\pi^{-(4m+2)}\Gamma(4m+2)L(4m+2,\left(\frac{d}{*}\right)).
\end{eqnarray}
As above we can use generalised Bernoulli numbers.
Hence by Hirzebruch-Mumford proportionality we obtain for $d>1$ that
$\dim S_k(\widetilde{\Orth}^+(N_{2d}^{(m)}))$ grows as 
$$
\frac{2^{-4m-1}}{(8m+2)!}
\frac{\,B_2\cdot B_4 \dots B_{8m+2}\,}{(8m+2)!!}\cdot
\frac{B_{4m+2,\chi_D}}{4m+2}\,
t^{8m+3}\prod_{p|t}(1-\chi_D(p)p^{-(4m+2)})\,k^{8m+2}.
$$
Here, as before, $d=d_0t^2$, with $d_0$ square-free, and $D=d_0$ is
the discriminant of the quadratic extension $\QQ(\sqrt{d})$.
For $d=1$ we have an extra factor $2$, $t=1$, $\chi_D\equiv 1$
and $B_{4m+2,\chi_D}=B_{4m+2}$.
In this case the lattice $N_{2d}^{(m)}$ is unimodular and the
formula again agrees with our previous computations in
Section \ref{sub_evenunimodular}.

\bibliographystyle{alpha}

\bigskip
\noindent
Valery Gritsenko\\
Universit\'e Lille 1\\
UFR de Math\'ematiques\\
F-59655 Villeneuve d'Ascq, Cedex\\
France\\
{\tt valery.gritsenko@math.univ-lille1.fr}
\bigskip

\noindent
K.~Hulek\\
Institut f\"ur Algebraische Geometrie\\
Leibniz Universit\"at Hannover\\
D-30060 Hannover\\
Germany\\
{\tt hulek@math.uni-hannover.de}
\bigskip

\noindent
G.K.~Sankaran\\
Department of Mathematical Sciences\\
University of Bath\\
Bath BA2 7AY\\
England\\
{\tt gks@maths.bath.ac.uk}

\end{document}